\documentclass{amsart} 

\newtheorem{theorem}{Theorem}[section]
\newtheorem*{maintheorem}{Theorem}
\newtheorem*{claim}{Claim}
\newtheorem{lemma}[theorem]{Lemma}
\newtheorem{proposition}[theorem]{Proposition}

\theoremstyle{definition}
\newtheorem{definition}[theorem]{Definition}
\newtheorem{remark}[theorem]{Remark}
\newtheorem{question}[theorem]{Question}
\newtheorem*{acknowledgments}{Acknowledgments}

\newcommand{\PP}{\mathbb{P}}

\newcommand{\GG}{\mathbb{G}}

\newcommand  {\shA}     {\mathcal{A}}
\newcommand  {\shB}     {\mathcal{B}}

\newcommand  {\shDiv}   {\mathcal{D} \!\text{\textit{iv}}}

\newcommand  {\shADiv}  {\mathcal{A}\mathcal{D}\!\text{\textit{iv}}}

\newcommand  {\shE}     {\mathcal{E}}
\newcommand  {\shF}     {\mathcal{F}}
\newcommand  {\shG}     {\mathcal{G}}

\newcommand  {\shAut}   {\mathcal{A}\!\text{\textit{ut}}}
\newcommand  {\shHom}   {\mathcal{H}\!\text{\textit{om}}}
\newcommand  {\shEnd}   {\mathcal{E}\!\text{\textit{nd}}}
\newcommand  {\shMat}   {\mathcal{M}\!\text{\textit{at}}}
\newcommand  {\shTor}   {\mathcal{T}\!\text{\textit{or}}}
\newcommand  {\shI}     {\mathcal{I}}
\newcommand  {\shJ}     {\mathcal{J}}

\newcommand  {\shM}     {\mathcal{M}}

\newcommand  {\shL}     {\mathcal{L}}

\newcommand  {\shS}     {\mathcal{S}}
\newcommand  {\shT}     {\mathcal{T}}

\newcommand  {\shP}     {\mathcal{P}}

\newcommand  {\shZ}     {\mathcal{Z}}

\newcommand  {\foY}     {\mathfrak{Y}}

\newcommand  {\Aut}     {\operatorname{Aut}}
\newcommand  {\ADiv}     {\operatorname{ADiv}}

\newcommand  {\Br}      {\operatorname{Br}}

\newcommand  {\Cl}      {\operatorname{Cl}}
\newcommand  {\cl}      {\operatorname{cl}}
\renewcommand{\cong}    {\equiv}

\newcommand  {\depth}   {\operatorname{depth}}

\newcommand  {\Div}     {\operatorname{Div}}

\newcommand  {\et}      {{\text{\'{e}t}}}

\newcommand  {\GL}      {\operatorname{GL}}
\newcommand  {\Grass}      {\operatorname{Grass}}

\newcommand  {\Isom}    {\operatorname{Isom}}
\renewcommand  {\k}     {\kappa}

\newcommand  {\liminv}  {\varprojlim}

\newcommand  {\lra}     {\longrightarrow}

\newcommand  {\NS}      {\operatorname{NS}}

\renewcommand{\O}       {\mathcal{O}}

\newcommand   {\hd}      {\operatorname{hd}}

\newcommand  {\Pic}     {\operatorname{Pic}}
\newcommand  {\PGL}     {\operatorname{PGL}}

\newcommand  {\ra}      {\rightarrow}

\newcommand  {\rank}    {\operatorname{rank}}

\newcommand  {\Sing}    {\operatorname{Sing}}

\newcommand  {\Spec}    {\operatorname{Spec}}

\newcommand  {\sh}      {\operatorname{sh}}

\def\mydate{\number\day\space\ifcase\month \or January\or February\or March\or 
April\or May\or\June\or\July\or
August\or September\or October\or November\or December\fi \space\number\year}

\begin{document}

\title[Azumaya algebras on surfaces]{There are enough Azumaya algebras on 
surfaces}

\author[Stefan Schroeer]{Stefan Schr\"oer}
\address{Mathematische Fakult\"at, Ruhr-Universit\"at, 
               44780 Bochum, Germany}
\curraddr{M.I.T. Mathematical Department, Room 2-155,
77 Massachusetts Avenue, Cambridge, MA 02139-4307,
USA}
\email{s.schroeer@ruhr-uni-bochum.de}



\dedicatory{Final version, 26 February 2001}

\begin{abstract}
Using Maruyama's theory of elementary transformations, I show   
that the Brauer 
group surjects onto the cohomological 
Brauer group  for separated geometrically
normal algebraic surfaces. As an application, I infer  the existence of   
nonfree 
vector bundles on proper 
normal algebraic surfaces.
\end{abstract}

\maketitle

\section*{Introduction}

Generalizing the classical theory of central simple algebras over fields, 
Grothen\-dieck \cite{GB} introduced the Brauer group $\Br(X)$ and the 
cohomological 
Brauer group $\Br'(X)$ for schemes. 

Let me recall the definitions.
The \emph{Brauer group} $\Br(X)$ comprises equivalence classes of Azumaya 
algebras.
Two  Azumaya algebras $\shA,\shB$ are called equivalent if there are 
everywhere nonzero 
vector bundles $\shE,\shF$ with $\shA\otimes\shEnd(\shE)\simeq 
\shB\otimes\shEnd(\shF)$.
Let us define the \emph{cohomological Brauer group} $\Br'(X)$ as the torsion 
part 
of 
the \'etale cohomology group $H^2(X,\GG_m)$. Nonabelian cohomology
gives  an inclusion 
$\Br(X)\subset \Br'(X)$, and Grothendieck  asked whether this is bijective.

It would be nice to know this for the following reason: The cohomological 
Brauer group is related to various other cohomology groups via exact 
sequences, 
and this 
is useful for computations. In contrast, it is almost impossible to calculate 
the 
Brauer group of a scheme directly from the definition. 
Here is a list of schemes with $\Br(X)=\Br'(X)$:
\begin{enumerate}
\item 
Schemes of dimension $\leq 1$ and regular surfaces (Grothendieck \cite{GB}).
\item 
Abelian varieties (Hoobler \cite{Hoobler 1972}).
\item 
The union of two affine schemes with affine intersection (Gabber \cite{Gabber 
1980}).
\item 
Smooth toric varieties (DeMeyer and Ford \cite{Demeyer; Ford 1993}).
\end{enumerate}

\noindent
On the other hand,   
a nonseparated normal surface   with $\Br(X)\neq\Br'(X)$ recently appeared in
\cite{Edidin et al 2000}. I wonder how the final answer to this puzzle  will 
look   
like. 
The goal of this paper is to
prove the following  Theorem.

\begin{maintheorem}
For separated geometrically normal algebraic surfaces, the inclusion
$\Br(X)\subset\Br'(X)$ is a bijection.
\end{maintheorem}

This adds some    singular and nonprojective schemes   to the preceding
list. For quasiprojective surfaces,
Hoobler (\cite{Hoobler 1982} Cor.\ 9)  deduced the result directly  from 
Gabber's 
Theorem on 
affine schemes. Without ample line bundles, a   different approach 
is required. Indeed, my initial motivation was to disprove the Theorem, rather 
than 
to prove it. The new idea is to use Maruyama's theory of elementary 
transformations. 

Here is an application of the preceding result:

\begin{maintheorem}
Each proper normal algebraic surface admits a nonfree vector bundle.
\end{maintheorem}

It might easily happen that all line bundles are free 
\cite{Schroeer 1999a}.
The existence of nonfree vector 
bundles  can  be viewed as a
generalization, in dimension  two, of  Winkelmann's Theorem \cite{Winkelmann 
1993}, which
asserts that each   compact 
complex manifold has   nonfree holomorphic vector bundles.

This paper has four sections. In the first section, I relate Azumaya algebras 
that are trivial on large open subsets to certain reflexive sheaves. 
 In Section 2, we turn to normal surfaces and  construct 
Azumaya
algebras that  are generically trivial by constructing the corresponding 
reflexive
sheaves. This prepares the proof of the main Theorem, which appears in Section 
3. The   
idea in the proof is to apply   elementary 
transformations to Brauer--Severi schemes.
The last section contains the  existence result for nonfree vector bundles. 

\begin{acknowledgments}
This research was done  in Bologna,
and I am grateful to the   Mathematical Department for its hospitality.
I wish to thank Angelo Vistoli for suggestions, encouragement, and many 
stimulating discussions. Furthermore, I
wish to thank Ofer Gabber, the referee, for his precise report.
He found and corrected  several mistakes. Several crucial steps
are entirely due to Gabber, and the paper would be impossible
without his contribution.
The revision was done at M.I.T., and I wish to thank the Mathematical
Department for its hospitality. Finally, I thank the DFG for financial
support.
\end{acknowledgments}

\section{Azumaya algebras via reflexive sheaves}
\label{Azumaya reflexive}

In this section, we shall describe Azumaya algebras that have
trivial Brauer class on certain large open subsets. Throughout, $X$ will be 
a noetherian scheme. Let us call an open subset $U\subset X$ \emph{thick}
if it contains all points  $x\in X$ with  $\depth(\O_{X,x})\leq 1$.
In other words,  
$\depth_{X-U}(\O_X)\geq 2$. A coherent $\O_X$-module $\shF$ is called 
\emph{almost locally free} if it is locally free on some thick
open subset $U\subset X$, and has $\depth_{X-U}(\shF)\geq 2$.
Such sheaves behave well under suitable 
restriction and extension functors:

\begin{lemma}
\label{equivalence}
Let $i:Y\subset X$ be a thick open subset. Then the restriction map 
$\shF\mapsto i^*(\shF)$
and the direct image map $\shG\mapsto i_*(\shG)$ induce inverse equivalences
between the  categories of almost locally free sheaves
on $X$ and $Y$, respectively.
\end{lemma}

\proof
This is similar to the 
proof of \cite{Hartshorne 1994} Theorem 1.12. Fix an almost locally
free $\O_X$-module $\shF$.
First, we check that $\Gamma(V,\shF)\ra\Gamma(V\cap Y,\shF)$ 
is bijective for all affine  open subsets $V\subset X$.
Setting $A=V-V\cap Y$, we have an exact sequence of local cohomology groups 
\begin{equation}
\label{local cohomology sequence}
0\lra H^0_A(V,\shF)\lra H^0(V,\shF)\lra 
H^0(V\cap Y,\shF)\lra H^1_A(V,\shF)\lra 0.
\end{equation}
Since $\depth_A(\shF)\geq 2$, 
the cohomology groups with supports vanish by
\cite{Hartshorne 1967} Theorem 3.8. 
Therefore, the map in the middle is bijective.
As a consequence, the adjunction map $\shF\ra i_*i^*(\shF)$ is bijective, 
so that
the restriction functor $\shF\mapsto i^*(\shF)$ is fully faithful.

Second, we check that the functor $\shF\mapsto i^*(\shF)$ is essentially
surjective. Fix an almost locally free $\O_Y$-module $\shG$.
By \cite{EGA I} Corollary 6.9.8, the sheaf $\shG$ extends to  
a coherent $\O_X$-module $\shM$.
I claim that $\shF=\shM^{\vee\vee}$ is almost locally free.
This is a local problem, so we may assume that there is a partial
resolution $\shL_1\ra\shL_0\ra\shM^\vee\ra 0$ with 
coherent locally free sheaves, hence an exact sequence
\begin{equation}
\label{resolution}
0\lra \shF\lra\shL_0^\vee\lra\shL_1^\vee.
\end{equation}
Let $A=X-U$, where  $U\subset Y$ is a thick open subset
on which $\shG$ is locally free.
The exact sequence (\ref{resolution}) gives 
an inclusion $H^0_A(X,\shL_0^\vee/\shF)\subset H^0_A(X,\shL_1^\vee)$
and an exact sequence of local cohomology groups
$$
H^0_A(X,\shL_0^\vee/\shF)\lra H^1_A(X,\shF)\lra  H^1_A(X,\shL_0^\vee).
$$
Since $\depth_A(\O_X)\geq 2$, the outer groups vanish, and we conclude 
$\depth_A(\shF)\geq 2$.
Consequently, $\shF$ is almost locally free.
By the same argument, we see that $\shG^{\vee\vee}$ is almost locally free.
Since the canonical map  $\shG\ra\shG^{\vee\vee}$ is bijective 
on some thick open subset, we conclude that
it is bijective on $Y$, hence  $\shG\simeq i^*(\shF)$.

It remains to check that $\shG\mapsto i_*(\shG)$ is the 
desired inverse equivalence.
Extend $\shG$ to an  almost locally free $\O_X$-module $\shF$. 
Since the  adjunction map $\shF\ra i_*i^*(\shF)$ is
bijective, we are done.
\qed

\begin{remark}
Given two almost locally free sheaves $\shF_1$ and $\shF_2$, 
the sheaf $\shF=\shHom(\shF_1,\shF_2)$
is almost locally free as well. This is because,
by Lemma \ref{equivalence}, the middle map  
in (\ref{local cohomology sequence}) is bijective, so that the 
cohomology groups with support
vanish. As a consequence, almost locally free sheaves are reflexive.
\end{remark}

An Azumaya algebra $\shA$ is called \emph{almost trivial} if its Brauer class 
$\cl(\shA)\in\Br(X)$ 
vanishes on some
thick open subset. This easily implies that 
$\shA\simeq\shEnd(\shF)$ for some almost locally free sheaf $\shF$.
However, the condition that the $\O_X$-algebra
$\shEnd(\shF)$ is an Azumaya algebra implies more.
 
\begin{definition}
A coherent $\O_X$-module $\shF$ is called \emph{balanced} if for each
geometric point 
$\bar{x}\ra X$,
there is a decomposition $\shF_{\bar{x}}\simeq
\bigoplus_{i=1}^r\shL_{\bar{x}}$ with $r>0$  for some 
almost invertible $\O_{X,\bar{x}}$-module 
$\shL_{\bar{x}}$.
\end{definition}

Here $\O_{X,\bar{x}}$ is the strict henselization of the local ring 
$\O_{X,x}$.
Perhaps it goes without saying that \emph{almost invertible sheaves} are 
invertible on a thick
open subset and have depth $\geq 2$ outside. By fpqc-descent, 
balanced
sheaves are almost locally free. They are closely related to Azumaya algebras,
and the following result reduces the existence of certain Azumaya algebras to
the existence of balanced sheaves.

\begin{proposition}
Let $\shF$ be a balanced $\O_X$-module. Then $\shEnd(\shF)$ is an almost 
trivial Azumaya algebra,
and each almost trivial Azumaya algebra has  this form.
\end{proposition}

\proof
Obviously, $\shA=\shEnd(\shF)$ is a trivial 
Azumaya algebra on some thick open subset.
To check that it is an Azumaya algebra on $X$, 
we may assume that $X$ is strictly local,
and that $\shF=\bigoplus_{i=1}^r\shL$ for some almost invertible sheaf $\shL$. 
Being bijective on some thick
open subset, the map $\O_X\ra\shEnd(\shL)$ is everywhere  bijective by 
Lemma \ref{equivalence}. Consequently,
$\shA\simeq\shMat_r(\O_X)$ is an Azumaya algebra.

Conversely, let $\shA$ be an almost trivial Azumaya algebra. 
Choose a thick open subset 
$i:U\subset X$ on which the Brauer class is trivial. 
Then there is an isomorphism 
$\shA_U\ra\shEnd(\shG)$
for some locally free $\O_U$-module $\shG$. By Lemma \ref{equivalence}, 
this induces
an isomorphism of algebras $\shA\ra\shEnd(\shF)$, where $\shF=i_*(\shG)$.

If remains to check that the almost locally free
sheaf $\shF$ is balanced. To do so, we may 
assume that $X$ is strictly local, so $\shA=\shMat_r(\O_X)$.
Now $\shF$ is a module over $\shMat_r(\O_X)$. By Morita equivalence 
(see e.g.\ \cite{Knus 1991} p.~53), 
$\shF=\bigoplus_{i=1}^r \shM$ for some coherent $\O_X$-module $\shM$. Clearly, 
$\shM$ is invertible on a thick open subset and has depth
$\geq 2$ outside. In other words, $\shM$ is almost invertible.
\qed

\medskip
Next, we shall generalize some notions from Hartshorne's paper on
generalized  divisors
\cite{Hartshorne 1994}. For simplicity, we 
assume that $X$ satisfies Serre's condition $(S_2)$, such
that the points of codimension one are precisely the points
of depth  one. Set 
$$
\shADiv_X = \bigoplus _{x\in X^{(1)}} (i_x)_*(\Div(\O_{X,x}))
$$
where the sum runs over all points of codimension one.
The elements of the group $\ADiv(X)=\Gamma(X,\shADiv_X)$ are called 
\emph{almost Cartier  divisors}.
For normal schemes, $\shADiv_X$ is just the sheaf of Weil divisors. As in the 
normal case, an almost Cartier 
divisor $D\in\ADiv(X)$ defines an almost invertible sheaf $\O_X(D)$, which is 
invertible in codimension one and satisfies Serre's condition $(S_2)$.
According to \cite{EGA IVd} Proposition 21.1.8, the canonical map 
$\shDiv_X\ra \shADiv_X$ is injective. The exact sequence
$$
0\lra \shDiv_X\lra \shADiv_X\lra\shP_X\lra 0
$$ 
defines an abelian sheaf $\shP_X$ on the \'etale site $X_\et$.
For a geometric point $\bar{x}\ra X$ with corresponding strict localization 
$\O_{X,\bar{x}}=\O^{\sh}_{X,x}$, the stalk is 
$$
\shP_{X,\bar{x}} = \ADiv(\O_{X,\bar{x}}) / \Div(\O_{X,\bar{x}}).
$$
For normal schemes, this reduces to the class group $\Cl(\O_{X,\bar{x}})$.
The preceding short exact sequence gives a long exact sequence in \'etale 
cohomology
$$
\ADiv(X)\ra H^0(X,\shP)\lra H^1(X,\shDiv_X)\lra H^1(X,\shADiv_X).
$$
We also have an exact sequence
$$
0\lra H^1(X,\shDiv_X)\lra H^2(X,\GG_m)\lra H^2(X^{(0)},\GG_m),
$$
where $X^{(0)}=\coprod\Spec(\O_{X,\eta})$ is the scheme of generic points, and
you easily 
infer that the classes of   almost trivial Azumaya algebra lie in the image
of the iterated coboundary map
\begin{equation}
\label{coboundary}
H^0(X,\shP)\lra H^1(X,\shDiv_X)\lra H^2(X,\GG_m).
\end{equation}
If $\shF$ is an almost locally free sheaf 
with decompositions 
$\shF_{\bar{x}}=\bigoplus_{i=1}^r\shL_{\bar{x}}$,   the function 
$$
\bar{x} \mapsto \cl(\shL_{\bar{x}})\in 
\ADiv(\O_{X,\bar{x}})/\Div(\O_{X,\bar{x}}) = 
\shP_{X,\bar{x}}
$$ 
depends only on $\shF$ and
yields a section 
$s_\shF\in\Gamma(X,\shP_X)$. The following result is
due to  Gabber:

\begin{proposition}[Gabber]
\label{section}
Let  $\shF$ be an balanced $\O_X$-module. Then the 
class  of the Azumaya algebra
$\shEnd(\shF)$ in $H^2(X,\GG_m)$ 
is the inverse of the image of the section $s_\shF\in H^0(X,\shP_X)$
under the iterated coboundary map in (\ref{coboundary}).
\end{proposition}

\proof
Set $\shA=\shEnd(\shF)$. 
According to \cite{Giraud 1971} p.\ 341, its  cohomology class in
$H^2(X,\GG_m)$ is given by the gerbe $d(\shA$) 
of  trivializations for $\shA$, which associates
to each \'etale $U\ra X$ the groupoid of  pairs $(\shE,\varphi)$, where 
$\shE$ is a locally free
$\O_U$-module, and $\varphi:\shEnd(\shE)\ra\shA_U$ is an isomorphism.
The action $\GG_m\ra\shAut(\shE,\varphi)$ is given by homotheties on $\shE$.

Set $\shM^\times_X=g_*g^*(\GG_m)$, where $g:X^{(0)}\ra X$ is the inclusion
of generic points,
 and let $\pi:\shM^\times_X\ra\shDiv_X$ and $p:\shADiv_X\ra\shP_X$ 
be the natural surjections. 
Then the image of $s_\shF\in H^0(X,\shP_X)$ under the iterated
coboundary map is given
by the gerbe of $\shM^\times_X$-liftings of the
$\shDiv_X$-torsor $p^{-1}(s_\shF)$.
This gerbe associates to each \'etale $U\ra X$ the groupoid of  pairs
$(\shT,\psi)$, where $\shT$ is a $\shM^\times_U$-torsor, and
$\psi:\shT\ra p^{-1}(s_\shF)$ is a $\pi$-morphism of torsors.
Moreover, the action $\GG_m\ra\shAut(\shT,\psi)$ is given 
by translation on $\shT$.

We have to construct an equivalence between the preceding two stacks
that is equivariant for the sign change  map $-1:\GG_m\ra\GG_m$. 
First note that the stack of pairs $(\shE,\varphi)$ is 
$\GG_m$-equivalent to the 
stack of triples $(\shE,\shL,\psi)$, where $\shE$ is a locally free 
$\O_U$-module,
and $\shL$ is an almost invertible $\O_U$-module, and 
$\psi:\shE\otimes\shL\ra\shF$ is an isomorphism.
The action $\GG_m\ra\shAut(\shE,\shL,\psi)$ is given by homotheties on 
$\shE$ and inverse homotheties
on $\shL$.

Obviously, the local  section $s_\shL\in\Gamma(U,\shP_X)$ 
is nothing but the restriction of the global section
$s_\shF\in\Gamma(X,\shP_X)$. Now let $\shM^\times_U(\shL)$ be the sheaf of 
invertible meromorphic sections
in $\shL|_U$, which is a $\shM^\times_U$-torsor. Dividing out the induced 
$\GG_m$-action,
we obtain a $\pi$-morphism 
$$
\varphi:\shM^\times_U(\shL)\lra p^{-1}(s_\shL)=p^{-1}(s_\shF|_U).
$$
Consequently, the functor
$(\shE,\shL,\psi)\mapsto (\shM^\times_U(\shL),\varphi)$
gives the desired anti\-equivariant equivalence of $\GG_m$-gerbes.
\qed

\medskip
For the rest of the section, we assume that $X$ satisfies Serre's condition
$(S_1)$, that is,  $X$ has no embedded components. Let 
$g:X^{(0)}\ra X$ be the
inclusion of the generic points. We can relate 
generically trivial Brauer classes to torsors of Cartier divisors as follows.
The exact sequence
$$
0\lra \GG_m\lra g_*g^* (\GG_m) \lra \shDiv_{X} \lra 0,
$$
together with $R^1g_*(\GG_{m,X^{(0)}})=0$ and $\Pic(X^{(0)})=0$, 
gives an exact sequence
\begin{equation}
\label{couboundary sequence}
0\lra  H^1(X,\shDiv_{X})\lra H^2(X,\GG_m)\lra H^2(X^{(0)},\GG_m).
\end{equation}
Hence each $\alpha\in\Br(X)$ with $g^*(\alpha)=0$ comes from a 
$\shDiv_{X}$-torsor. To make this explicit,
choose an Azumaya algebra $\shA$, say of rank $r$,
 representing the class $\alpha$, and let 
$f:B\ra X$ be the associated \emph{Brauer--Severi} scheme. This is the
$\PP^{r-1}$-bundle 
$$
B=\Isom(\shMat_r(\O_X),\shA)\times_{\PGL_r}\PP^{r-1} 
$$ 
on the \'etale site $X_{\et}$. Here  we use the left $\PGL_r$-action
coming from the canonical representation 
$\PGL_r\ra \Aut(\PP^{r-1})$ described in 
\cite{Mumford; Fogarty; Kirwan 1993} Chap.\ 0 \S 5.
Set $P=B\times_X X^{(0)}$, and pick an invertible $\O_P$-module $\O_P(1)$ 
of fiber degree one.
This leads to a $\shDiv_X$-torsor $\shT$ as follows. Define
$$
\Gamma(X,\shT) = \left\{\cl(\shL,t)\right\},
$$
where $\shL$ is an invertible $\O_{B}$-module, and 
$t:\O_{P}(1)\ra\shL|_{P}$ is an isomorphism. Here 
$\cl(\shL,t)$ 
denotes  isomorphism class, and two
pairs $(\shL,t)$ and $(\shL',t')$ 
are called isomorphic if there is an isomorphism 
$\phi:\shL\ra\shL'$ with $t\circ \phi=t'$. 
Define $\shT $ in the same way on the \'etale site.

By \cite{EGA IVd} Proposition 21.2.11, 
the sections of $\shDiv_{X}$ correspond to $\cl(\shM,s)$, 
where $\shM$ is an invertible
$\O_X$-module, and $s:\O_{X^{(0)}}\ra\shM|_{X^{(0)}}$ is a trivialization. 
The map 
$$
(\shM,s),(\shL,t)\mapsto (f^*(\shM)\otimes\shL,f^*(s)\otimes t)
$$
turns $\shT$ into a $\shDiv_X$-torsor.

\begin{proposition}
The Brauer class $\cl(\shA)\in H^2(X,\GG_m)$ is the
opposite for the image of the 
torsor class $\cl(\shT)\in H^1(X,\shDiv_{X})$
under the coboundary map in (\ref{couboundary sequence}).
\end{proposition}

\proof
As in the proof of Proposition \ref{section},
the Brauer class is given by the gerbe of trivializations 
$\varphi:\shEnd(\shE)\ra \shA$. This gerbe is equivalent
to the gerbe of trivializations $(\shE,u)$, where
$u:B\ra \PP(\shE^\vee)$ is an isomorphism. 
According to \cite{Giraud 1971} Chap.\ V Lemma 4.8.1, 
the latter  gerbe is antiequivalent to the 
$\GG_m$-gerbe of invertible $\O_P$-modules $\shL$ of fiber degree one.

Set $\shM_X^\times=g_*g^*\GG_m$. The image of the torsor class
$\cl(\shT)\in H^1(X,\shDiv_X)$ is the gerbe of 
$\shM_{X}^\times$-liftings $\psi:\shS \ra \shT$ for 
the $\shDiv_{X}$-torsor $\shT$.
Given an invertible $\O_P$-module $\shL$ of fiber degree one,
we obtain an $\shM_X^\times$-torsor 
$$
\shS = \left\{(\shL,t) \mid t:\O_P(1)\stackrel{\simeq}{\lra}\shL_P\right\},
$$
where the  action is by multiplication on $t$.
Clearly, the quotient 
$\shS/\GG_m= \left\{\cl(\shL,t)\right\}$ is canonically isomorphic to $\shT$, 
so we obtain
a morphism of torsors $\psi:\shS\ra\shT$.
To see that $\shL\mapsto(\shS,\psi)$ is a $\GG_m$-equivalence, note that 
a $\shM_X^\times$-lifting of the torsor $\shT$ exists if and only if $\shT$ 
is trivial, because 
$H^1(X,\shM_X^\times)=0$.
\qed

\section{Generically trivial Brauer classes}

In this section, we turn to normal
 surfaces. The task is to   prove the following 
result, which is a major step towards showing $\Br(X)=\Br'(X)$.

\begin{proposition}\label{trivial}
Let $X$ be a separated normal algebraic surface.
Then $\Br(X)$ contains each class $\alpha\in\Br'(X)$ that is generically 
trivial.
\end{proposition}

\proof
We start with some preliminary reductions.
By \cite{Gabber 1980} Chap.\ II Lemma 4, we may assume that
the ground field $k$ is separably closed. 
Since $X$ is separated
and of finite type, the Nagata Compactification Theorem \cite{Nagata 1962}
gives a compactification $X\subset \bar{X}$.
By resolution of singularities, we may assume that $\Sing(X)=\Sing(\bar{X})$.
As in \cite{GB} Chap.\ II Theorem 2.1, each generically 
trivial Brauer class $\alpha\in\Br'(X)$
extends to $\Br'(\bar{X})$, so we may begin the proof  with the additional
assumption that $X$ is proper.

As discussed in Section \ref{Azumaya reflexive}, the exact 
sequence 
$$
0\lra H^1(X,\shDiv_X)\ra H^2(X,\GG_m) \lra H^2(X^{(0)},\GG_m)
$$
shows that our cohomology class $\alpha\in\Br'(X)$ lies in $H^1(X,\shDiv_X)$.
The exact sequence 
$$
0\lra \shDiv_X\lra \shZ^1_X\lra \shP_X\lra 0
$$
yields an exact sequence
$$
Z^1(X)\lra H^0(X,\shP_X)\lra H^1(X,\shDiv_X)\lra 0.
$$
Choose a global section $s\in H^0(X,\shP_X)$ mapping to $\alpha$. Since  
$r\alpha=0$ 
for some
integer $r> 0$, there is a global Weil  divisor $E\in Z^1(X)$ mapping to the 
section
$rs\in H^0(X,\shP_X)$. 

The sheaf $\shP_X$ is supported by the singular locus $\Sing(X)$. For each 
singular point
$x\in X$,  the stalk is $\shP_{X,x}=\Cl(\O_{X,x}^h)$, where $\O_{X,x}\subset 
\O_{X,x}^h$
is  the henselization (which is  the strict localization, because $k$ is 
separably 
closed). Suppose $x_1,\ldots, x_m\in X$ are the singularities, and set 
$$
X^h=\Spec(\prod_{i=1}^m\O_{X,x_i}^h) = \coprod_{i=1}^m\Spec(\O_{X,x_i}^h).
$$
Then $H^0(X,\shP_X)=\Cl(X^h)$. Choose a Weil divisor $D\in Z^1(X^h)$ 
representing the section $s\in H^0(X,\shP_X)$, such that
$E|_{X^h}\sim rD$.
According to Proposition \ref{section}, it suffices to 
construct a reflexive 
$\O_X$-module $\shF$ with
$$
\shF\otimes\O_{X^h}=\bigoplus_{i=1}^r \O_{X^h}(-D),
$$
for then $\shA=\shEnd(\shF)$ would be the desired Azumaya algebra.

Let $f:Y\ra X$ be a resolution of singularities, and 
$Y_0\subset Y$ be the reduced exceptional curve.
Set $Y^h=Y\times_X X^h$. 
The following crucial result, which is due
to Gabber, tells us that certain vector bundles on $Y^h$
are already determined on suitable infinitesimal neighborhoods of $Y_0$.

\begin{lemma}[Gabber]
\label{exceptional}
Let $\shB$ be a family of locally free $\O_{Y_0}$-modules
of fixed rank $n\geq 0$. Suppose that $\shB$ 
is, up to tensoring with line bundles, a  bounded family.
Then there is an exceptional curve $R\subset Y$ so that
the  map $H^1(Y^h,\GL_n)\ra H^1(R,\GL_n)$ is bijective on the subsets
of vector bundles whose restriction to $Y_0$ lies in $\shB$.
\end{lemma}

\proof
Let $\shI\subset \O_Y$ be the ideal of $Y_0\subset Y$.
Since the intersection matrix for the irreducible
components of $Y_0$ is negative definite, there is an exceptional curve
$A\subset Y$ with support $Y_0$ so that $\O_{Y_0}(-A)$ is ample.
Then $\O_{A}(-A)$ is ample as well. Since the
family  $\left\{\shEnd(\shE)\mid \shE\in\shB\right\}$ 
is bounded,
there is an integer $m_0> 0$ so that
\begin{equation}
\label{vanishing}
H^1(A,\shEnd(\shE)\otimes\O_{A}(-mA))=0
\end{equation}
for all  $m\geq m_0$ and all $\shE\in\shB$.

Let $\foY\subset Y$ be the formal completion along $Y_0\subset Y$. 
We first check the statement of the Lemma
for the formal scheme $\foY$ instead of $Y^h$. 
Note that the  canonical map 
$$
H^1(\foY,\GL_n)\lra \liminv H^1(mA,\GL_n)
$$
is bijective, as explained in \cite{Artin 1969}, proof of Theorem 3.5.
Let $\shJ\subset \O_Y$ be the ideal of $A\subset Y$. 
The obstruction to lifting a
vector bundle $\shE$ on $mA$ to $(m+1)A$ lies in 
$$
H^2(R,\shEnd(\shE)\otimes \shJ^m/\shJ^{m+1})=0,
$$
so the restriction maps
$H^1(\foY,\GL_n)\ra  H^1(mA,\GL_n)$ are  surjective for all $m\geq 0$.

Fix an integer $m\geq m_0$, and
let $\shE,\shE'$ be two vector bundles on $(m+1)A$ that
are isomorphic on $mA$ and whose restrictions to $Y_0$ 
belong to the family $\shB$.
Choose an isomorphism $\psi:\shE|_{mA}\ra\shE'|_{mA}$. 
Locally on $(m+1)A$, 
we can lift this
isomorphism to an isomorphism $\shE\ra\shE'$. 
The sheaf of such liftings is a torsor
under 
$$
\shHom(\shE,\shE')\otimes\shJ^{m}/\shJ^{m+1}\simeq 
\shEnd(\shE)\otimes\shJ^{m}/\shJ^{m+1}\simeq
\shEnd(\shE)\otimes\O_A(-mA).
$$
This sheaf has no first cohomology by (\ref{vanishing}).
The upshot is that  a global 
lifting of  $\psi:\shE|_{mA}\ra\shE'|_{mA}$ exists. 
Consequently, for $R=m_0A$,  the  mapping 
$H^1(\foY,\GL_n)\ra H^1(R,\GL_n)$ is  bijective on 
the subsets of vector bundles  whose
restriction belongs to the family  $\shB$.

Finally, we pass  to $Y^h$.
By the Artin Approximation Theorem (\cite{Artin 1969} Thm.\ 3.5),
the map $H^1(Y^h,\GL_n)\ra 
H^1(\foY,\GL_n)$ is injective and has dense image.
So given a formal vector bundle $\shE$ with  $\shE|_{Y_0}\in\shB$,
we find a vector bundle $\shE^h$ on $Y^h$ with $\shE^h|_R\simeq\shE|_R$.
By the choice of $R$, this implies $\shE^h|_\foY\simeq \shE$.
\qed

\medskip
We proceed with the proof of   Proposition \ref{trivial}. Let
$\shB$ be the family of vector bundles on $Y_0$ of rank $r$ that are free
up  to tensoring with line bundles, and choose an exceptional divisor 
$R\subset Y$ as in the preceding Lemma.
 Let $D' \in Z^1(Y^h)$ be the
strict  transform of $D\in Z^1(X^h)$.
Let $E'\in Z^1(Y)$ be the unique Weil divisor that is the strict transform
of $E\in Z^1(X)$ on $Y-Y_0$, and satisfies $E'\sim rD'$ on $Y^h$.
Set 
$\shL
=\O_R(D')$. According to Lemma \ref{exceptional}, we have to construct a 
locally 
free 
$\O_Y$-module $\shE$ with $\shE_R\simeq \bigoplus _{i=1}^r\shL^\vee$. 
For then the 
double dual $\shF=f_*(\shE)^{\vee\vee}$ would be the desired reflexive 
$\O_X$-module.
We shall construct such a vector bundle as an elementary 
transformation of the 
trivial   bundle $\O_Y^{\oplus r}$.

Choose an ample divisor $A\subset Y$.
Replacing the divisors $D$ and $E$ by $D+ f_*(tA)$ and $E+f_*(rtA)$, 
respectively,  does not change
the class
$\alpha\in \Br'(X)$. Choosing 
$t\gg 0$, we may assume that 
$\shL=\O_R(D')$ is very ample, and that 
$H^1(Y,\O_Y(E'-R))=0$. 
Next, choose pairwise disjoint effective Cartier divisors 
$D_1,\dots, D_r\subset R$, each one representing  
$\shL$. Let 
$s_i\in\Gamma(R,\shL)$ be the corresponding sections. Regard their product
$s_1\otimes\ldots \otimes s_r$ as a section of 
$\O_R(E')$. By construction, the group on the right in the exact sequence 
$$
  H^0(Y,\O_Y(E'))   \lra  H^0(R,\O_R(E')) \lra H^1(Y,\O_Y(E'-R))
$$
is zero. Consequently, 
$E'\in Z^1(Y)$ is linearly equivalent to an effective divisor 
$H\subset Y$ with 
$H\cap R= D_1\cup\ldots\cup D_r$. 
Now choose a closed subset  $S\subset H-R$
so that
each Cartier divisor $D_i\subset H$ is principal on
$H-S$. For each $1\leq i\leq r$,
choose an exact sequence
\begin{equation}
\label{divisors}
0\lra \O_H(C_i)\stackrel{t_i}{\lra} \O_H\lra 
\bigoplus_{j\neq i} \O_{D_j}  \lra 0
\end{equation}
for certain Cartier divisors 
$C_i\in\Div(H)$ supported by $S$.
As explained in \cite{Milne 1980}, p.~152, it suffices
to construct the desired Azumaya $\O_X$-algebra on 
$X-f(S)$.
Hence it suffices to construct the desired locally free
$\O_Y$-module $\shE$ on $Y-S$,
and we may replace $X$, $Y$ by the complements $X-f(S)$, $Y-S$,
respectively.
Now the exact sequence (\ref{divisors}) induces an exact sequence
$$
 \O_Y\stackrel{t_i}{\lra}
 \O_H\lra \bigoplus_{j\neq i} \O_{D_j}  \lra 0.
$$
The map 
$t=(t_1,\dots,t_r): \bigoplus_{i=1}^r\O_Y \lra \O_H$ is surjective, 
because the   
$D_i$ are pairwise disjoint. 
The exact sequence 
$$
0 \lra \shE  \lra \bigoplus_{i=1}^r\O_Y \stackrel{t}{\lra} \O_H  \lra 0
$$
defines a locally free 
$\O_{Y}$-module 
$\shE$, because the cokernel
$\O_H$ has homological dimension 
$\hd(\O_H)=1$. Restricting to the curve 
$R\subset Y$, we obtain an exact sequence 
$$
\shTor^1_{\O_{Y}}(\O_H,\O_R)  \lra \shE_R  \lra \bigoplus_{i=1}^r\O_R  
\stackrel{t_R}{\lra}
\O_{H\cap R}  \lra 0.
$$
The   term on the left is zero, because the curves 
$H,R\subset Y$ have no common components. By construction, 
$t_i|_{D_j} = 0$  for $i\neq j$, so the induced surjection 
$t_R$ is   a diagonal matrix of the form
$$
t_R = 
\begin{pmatrix} t_1|_{D_1} & & 0 \\ & \ddots  &\\ 0 && t_r|_{D_r} 
\end{pmatrix}:
\bigoplus_{i=1}^r \O_R \lra \bigoplus_{i=1}^r \O_{D_i} = \O_{H\cap R}.
$$
Consequently 
$\shE_R\simeq\bigoplus_{i=1}^r \shL^\vee$. Hence $\shE$ is the  
desired 
locally
free 
$\O_Y$-module.
\qed

\begin{remark}
In my first proof of Proposition \ref{trivial},  I used
a result of Treger (\cite{Treger 1978} Prop.\ 3.5),
which states that the map $H^1(Y^h,\GL_n)\ra 
H^1(R,\GL_n)$ is bijective for a suitable
exceptional curve $R\subset Y$.
As Gabber pointed out, this statement 
is wrong for $n\geq 2$. His counterexample
goes as follows.

Let $X=\Spec(A)$ be a complete normal local surface
singularity, and $f:Y\ra X$ a resolution of singularity, 
and $Y_m\subset Y$ the
infinitesimal neighborhoods of the
exceptional curve $Y_0$. Suppose $n=2$ for simplicity,
and assume there is an exceptional curve $R$ as above.
Fix an ample invertible $\O_Y$-module $\shL$, set
$\shL_m=\shL|_{Y_m}$, and choose
$m\geq 0$ with $R\subset Y_m$. The exact sequence
$$
0\lra \shL^{-t}_{0}(-Y_m) \lra \shL^{-t}_{m+1} \lra \shL^{-t}_{m} \lra 0
$$
yields an  exact sequence
$$
H^0(Y_m,\shL^{-t}_{m})\ra 
H^1(Y_{0},\shL^{-t}_{0}(-Y_m))\ra 
H^1(Y_{m+1},\shL^{-t}_{m+1})\ra
H^1(Y_{m},\shL^{-t}_{m})\ra 0.
$$
For $t\gg 0$, the group on the left is zero, and 
$H^1(Y_{0},\shL^{-t}_0(-Y_m))$ is nonzero.
It follows that 
there is a nonzero class $\zeta\in H^1(Y,\shL^{-t})$ 
restricting to zero in $H^1(Y_{m},\shL^{-t}_{m})$.
This defines a nonsplit extension 
\begin{equation}
\label{extension}
0\lra \shL^{-t}\lra \shE\stackrel{\psi}{\lra}\O_Y\lra 0,
\end{equation}
which splits on $R$. By the defining property of the curve
$R$, there is a bijection
$\phi:\O_Y\oplus \shL^{-t}\ra\shE$.  The composition
$$
\O_Y\stackrel{\phi}{\lra}\shE\stackrel{\psi}{\lra}\O_Y
$$
is surjective on $Y_m$, because 
$H^0(Y_m,\shL_{m}^{-t})=0$. By the Nakayama Lemma, the
composition is surjective on the formal completion, and
hence on $Y$ as well.  So the extension
(\ref{extension}) splits, contradicting  $\zeta\neq 0$.
\end{remark}

\section{Elementary transformations of Brauer--Severi schemes}

We come to the main result of this paper.

\begin{theorem}
\label{equal}
Let $X$ be a separated geometrically normal algebraic surface.  
Then we have $\Br(X)=\Br'(X)$.
\end{theorem}

\proof
According to \cite{Gabber 1980} Chap.\ II Lemma 4, we may
assume that the ground field is algebraically closed.
Fix a 
class $\alpha\in\Br'(X)$. In light of Proposition \ref{trivial}, 
it suffices to 
construct an Azumaya algebra representing $\alpha\in\Br'(X)$ generically.
Choose a resolution of singularities $f:Y\ra X$, and let $\foY\subset Y$ 
be the formal completion along the reduced exceptional curve $Y_0\subset Y$. 
According to \cite{GB} Theorem 2.1,
there is an Azumaya $\O_Y$-algebra $\shA$ representing $f^*(\alpha)$.
The task now is to choose such an 
Azumaya algebra so that the formal vector bundle
$\shA|_\foY$ is trivial. For then, as explained in \cite{Milne 1980} p.\ 152,
the $\O_X$-algebra  $f_*(\shA)$ is an Azumaya algebra, which 
represents $\alpha\in\Br'(X)$ generically.
Note that  we may remove finitely many closed smooth points
from $X$.

First, we check that the formal Azumaya algebra 
$\shA|_\foY$ is trivial. 
Since the ground field is separably closed, there is a locally
free $\O_{Y_0}$-module $\shE_0$ and an isomorphism 
$\varphi_0:\shEnd(\shE_0)\ra\shA|_{Y_0}$.
The following argument due to Gabber shows that 
the pair $(\shE_0,\varphi_0)$ extends over all infinitesimal neighborhoods 
$Y_0\subset Y_n$.
Let $d(\shA|_{Y_n})$ be the $\GG_m$-gerbe of trivializations of 
$\shA|_{Y_n}$.
The restriction map gives a cartesian functor 
$d(\shA|_{Y_{n+1}})\ra d(\shA|_{Y_n})$.
Consider the corresponding stack of liftings of trivializations.
This is a gerbe for the abelian sheaf $\shI^{n+1}/\shI^{n+2}$ 
on the \'etale site
of $Y_0$, where $\shI=\O_Y(-Y_0)$. Since $H^2(Y_0,\shI^{n+1}/\shI^{n+2})=0$, 
the gerbe of liftings is trivial.
Consequently, we have $\shA|_\foY\simeq\shEnd(\shE)$ for some locally 
free $\O_\foY$-module
 $\shE$.

The idea now is to make an elementary transformation along a curve
$H\subset Y$, so that $\shE$ becomes free on certain infinitesimal
neighborhood $Y_0\subset R$. Furthermore, we shall choose the curve
$R$ so that  the freeness of $\shE_R$ implies the freeness
of $\shE$.
This requires some preparation. 
Let me introduce  three numbers
$m,k,q$ depending on $Y$ and $\shE$.
First, set $m=\rank(\shE)$. 
Second, let $k\geq 1$ be the order of the cokernel for 
the map 
$\Pic(Y)\ra \NS(Y_0)$ onto the  N\'eron--Severi group. 
Third,  define
$q=1$ in characteristic zero.
In characteristic $p>0$, let $q>0$ be a $p$-th power so that the unipotent
part of $\Pic^0(\foY)$ is $q$-torsion. This works, because $\Pic^0_{\foY}$
 is an algebraic
group scheme. 

Now set  $r=mkq$, and  let $\shB$ be the family
of locally free $\O_{Y_0}$-modules of rank $r$ which are free
 up to tensoring with line bundles. Choose a curve $R\subset Y$
as in Lemma \ref{exceptional}. 
Finally, let $X\subset \bar{X}$ be a compactification with
$\Sing(X)=\Sing(\bar{X})$,
and let $Y\subset \bar{Y}$ be the 
corresponding compactification.

\begin{claim}
We can modify the Azumaya algebra $\shA$ and the vector bundle $\shE$ so
 that $\shE$ is a globally generated formal vector bundle of rank $r$, 
and that there is a very ample 
invertible
$\O_{\bar{Y}}$-module $\shL$ with $\det(\shE) = \shL_\foY $  and 
$H^1(\bar{Y},\shL(-R))=0$.
\end{claim}

\proof
Tensoring $\shE$ with an ample line bundle, we 
 archive that $\shE$ is globally generated and 
that $\det(\shE)$ is ample.
Next, we replace the Azumaya algebra 
$\shA$ by $\shA\otimes\shEnd(\O_Y^{\oplus k})$.  
This replaces the vector bundle $\shE$ by $\shE^{\oplus k}$, and $\det(\shE)$ 
by $ \det(\shE) ^{\otimes k}$.
Consequently, we may assume that $\det(\shE)$ is numerically equivalent to 
some ample invertible $\O_{\foY}$-module $\shL_\foY$.
Twisting $\shE$ by a suitable power of $\shL_\foY$,
we may assume that $\shL_\foY$ extends to an ample invertible 
$\O_{\bar{Y}}$-module $\shL$ with
$H^1(\bar{Y},\shL^{\otimes s}(-R))=0$ for all integers
$s>0$.

In characteristic zero, 
$\Pic^0(\foY)$ is a divisible group. Hence we find an invertible 
$\O_\foY$-module 
$\shM$ with 
$\shM^{\otimes mk} \otimes \det(\shE)\simeq \shL_\foY$. Replacing 
$\shE$ by 
$\shE\otimes\shM$, we have 
$\det(\shE)=\shL_\foY$. 
Now assume that we are in characteristic 
$p>0$. Set 
$G=\Pic^0_{\foY/k}$, and let 
$G'\subset G$ be the unipotent part. Then the quotient 
$G''=G/G'$ is semiabelian, and 
$G''(k)$ is a divisible group. As in characteristic zero, we may assume that 
$\det(\shE)\otimes\shL^\vee_\foY$ lies in the unipotent part of 
$\Pic^0(\foY)$, which is a $q$-torsion group. Passing to 
$\shA\otimes\shEnd(\O_X^{\oplus q})$ and 
$\shE^{\oplus q}$, we are done. 
\qed

\medskip
We continue proving  Theorem \ref{equal}. 
Set $r=\rank(\shE)$ and $\Gamma=\Gamma(R,\shE_R)$.   The canonical surjection 
$\Gamma\otimes\O_R\ra \shE_R$ yields a morphism 
$\varphi:R\ra\Grass_r(\Gamma)$ into the 
Grassmannian of r-dimensional quotients. 
Choose a
generic 
$r$-dimensional  subvector space $\Gamma'\subset \Gamma$.
For each integer $k\geq 0$, let 
$G_k\subset\Grass_r(\Gamma)$ be the subscheme of surjections 
$\Gamma\ra\Gamma''$ 
such
that  the composition
$\Gamma'\ra\Gamma''$ has rank 
$\leq k$. Note that 
$G_{r-1}$ is a reduced Cartier divisor, and that 
$G_{r-2}$ has codimension four (see \cite{Arbarello et al. 1985}, Sec.\ II.2).

By the dimensional part of Kleiman's Transversality Theorem
(\cite{Kleiman 1974} Thm.\  2), 
which is valid in all characteristics,  the map 
$\varphi:R\ra\Grass_r(\Gamma)$ is disjoint to 
$G_{r-2}$ and passes through 
$G_{r-1}$ in finitely many points.
The upshot of   this is that the quotient of the canonical map
$ \Gamma'\otimes \O_X\ra \shE$ is an invertible sheaf on some Cartier divisor 
$D\subset R$. Consequently, we have constructed an exact sequence 
$$
0 \lra \O_R^{\oplus r}  \lra \shE_R  \lra \O_D  \lra 0.
$$
In other words, the trivial vector bundle is the elementary transformation of 
$\shE$ with respect to the surjection 
$\shE_R  \ra \O_D$. In geometric terms: 
Blowing up 
$\PP(\shE_R)$ along the section 
$\PP(\O_D)\subset\PP(\shE_D)$ and contracting the strict transform of 
$\PP(\shE_D)$ yields 
$\PP^{r-1}_R$ (see \cite{Maruyama 1981} Thm.\ 1.4).

We  seek to extend this elementary transformation from the curve to the 
surface. Note that 
$\shL_R=\det(\shE_R)=\O_R(D)$. The exact sequence
$$
H^0(\bar{Y},\shL)\lra H^0(R,\shL_R)\lra H^1(\bar{Y},\shL(-R))=0
$$ 
implies that   
$D=H\cap R$ for some ample curve 
$H\subset Y$. Since the ground field 
$k$ is algebraically closed, Tsen's Theorem   gives 
$H^2(H,\GG_m)=0$. Removing finitely many smooth points from the open subset
$X\subset \bar{X}$,
we may assume  that $\shA_H=\shEnd(\O_H^{\oplus r})$.
In other words, if $P\ra Y$ is the Brauer--Severi scheme corresponding to
$\shA$, we have $P_H=\PP_H^{r-1}$.
Let $A$ be the semilocal ring of the 
curve $H$ corresponding to the closed points $D\subset H$.
The section $\PP(\O_D)\subset P_D$ is given by a surjection
$A^{\oplus r}\ra A/I$, where $I\subset A$ is the ideal of $D$.
By  Nakayama's Lemma,  this lifts to a surjection $A^{\oplus r}\ra A$.
So, if we shrink $X$ further, we can extend the
section 
$\PP(\O_D)\subset P_D$ to a section $S\subset P_H$.

Let $h:\hat{P}\ra P$ be the blowing-up with center $S\subset P$, 
and let $E\subset \hat{P}$ be the strict transform of the Cartier divisor 
$P_H\subset P$. I claim that there is a birational contraction 
$\hat{P}\ra P'$ contracting precisely the fibers of $E\ra H$ to points 
such that $P'$ is a Brauer--Severi scheme.  Over suitable \'etale 
neighborhoods, this follows from \cite{Maruyama 1981} Theorem 1.4. You easily 
check that 
these contractions glue together and define a contraction in the category of 
schemes.

By construction, the new Brauer--Severi scheme $P'\ra X$ has a trivial 
restriction $P'_R=\PP^{r-1}_R$. 
If $\shA'$ is the Azumaya algebra corresponding to the Brauer--Severi scheme 
$P'$, this implies $\shA'_R=\shEnd(\O^{\oplus r}_R)$. By the 
choice of the curve $R\subset Y$, this forces the formal vector bundle
$\shA'_{\foY}$ to be free. 
Consequently, the direct image $f_*(\shA')$ is an Azumaya 
$\O_{X}$-algebra  representing the   class 
$\alpha\in\Br'(X)$ generically.
\qed

\begin{remark}
The proof   works for separated geometrically normal 2-dimen\-sional 
algebraic spaces as well. This is  because
their resolutions are schemes.
\end{remark}

\begin{question}
The  hypothesis of \emph{geometric} normality  annoys me.
What happens for separated normal 2-dimensional noetherian schemes that are of 
finite type 
over   nonperfect   fields, or   over Dedekind rings, 
or have no base ring at all?
\end{question}

\section{Existence of vector bundles}

Given a scheme  $X$, one might ask whether $X$ admits   a nonfree vector
bundle. In dimension two, we can use   Brauer groups   to obtain a
positive answer:

\begin{theorem}
 Let $X$ be a proper normal  surface over a field $k$. Then there is 
a locally free  $\O_X$-module of finite rank that is not free.
\end{theorem}

\proof
Seeking a contradiction, we assume that each vector bundle is free.
We may assume that $k=\Gamma(X,\O_X)$.
First, I reduce to the case that the ground field $k$ is separably closed.
Let $k\subset L$ be a separable closure, and let $\shE_L$ be a vector bundle 
on $X_L=X\otimes L$, say of rank $r\geq 0$. Then there is a finite separable 
field extension $k\subset K$, say of degree $d\geq 1$, such that $\shE_L$ 
comes from a vector bundle $\shE_K$ on $X_K$.
Let $p:X_K\ra X$ be the canonical projection. Then $\shF=p_*(\shE_K)$ is a 
vector bundle of rank $dr$, hence free by assumption. This gives 
$$
\Gamma(X_K,\shE_K) = \Gamma(X,\shF) \simeq k^{\oplus rd} \simeq K^{\oplus r}.
$$
Now you easily choose $r$ sections of $\shF$ that are linearly independent 
over $K$, which gives the desired trivialization of $\shE_K$.

From now on, assume that $k$ is separably closed.
Note that $\Pic(X)=0$, so $X$ is nonprojective, hence it must contain some 
singularities. Let $f:Y\ra X$ be  a resolution of singularities. Choose an 
exceptional divisor $R\subset X$ so that 
$\Pic(R)=\Pic(\foY)$, where $\foY\subset Y$ is the formal completion along 
the exceptional curve.  The spectral 
sequence for $\GG_{m,X}=f_*(\GG_{m,Y})$  gives an 
exact sequence 
$$
0\lra \Pic(Y)\lra \Pic(R)\lra H^2(X,\GG_m)\ra H^2(Y,\GG_m).
$$
Set $G=\Pic^0(R)/\Pic^0(Y)$, and let $H\subset \NS(Y)$ be the kernel of the 
restriction map   $\NS(Y)\ra \NS(R)$ for the N\'eron--Severi groups. 
The snake lemma gives an inclusion $G/H\subset H^2(X,\GG_m)$.
To proceed, we need a well-known fact:

\begin{lemma}
For each   $l>0$ prime to $\operatorname{char}(k)$, 
the group $\Pic^0(R)$ is $l$-divisible.
\end{lemma}

\proof
We have an exact
sequence
$$
0\lra\Pic^0(R)\lra\Pic^0_{R/k}(k)\lra \Br(k).
$$
Since $\Br(k)=0$, we have to see
that the multiplication 
morphism
$l:P\ra P$ is surjective on the smooth 
algebraic group scheme $P=\Pic^0_{R/k}$. 
Since $P$ is connected, it suffices to check that $l:P\ra P$ is 
open. The completion at the origin $0\in P$ is a formal group,
given by a formal power series ring $k[[X_1,\ldots,X_n]]$ together
with $n$ formal power series 
$$
F_i(X_1,\ldots,X_n,Y_1,\ldots,Y_n) = X_i+Y_i + 
\text{terms of higher order.}
$$
Multiplication by $l$
is given by $l-1$ substitutions
$$
[l]^*X_i = F_i([l-1]^*X_1,\ldots,[l-1]^*X_n,X_1,\ldots,X_n) \cong lX_i,
$$
modulo terms of higher order. Since $l$ is prime to the characteristic of the 
ground field,
this is bijective.
 Consequently, $l$ is \'etale on $\O_{P,0}^\wedge$,
hence \'etale on $P$, and therefore open.
\qed

\medskip
We continue with the proof of the Proposition. If $G=0$, then $X$ would have 
nontrivial line bundles (\cite{Schroeer 1999b}, Prop.\ 4.2), which is absurd.  
So $G$ is a 
nonzero $l$-divisible group. On the other hand,   $H $ is finitely generated. 
We conclude that   $G/H\subset H^2(X,\GG_m)$ contains many torsion 
points. Consequently, $\Br'(X)$   contains nonzero generically trivial 
classes. By Proposition
\ref{trivial}  there is a nontrivial Azumaya $\O_X$-algebra
$\shA$.

Setting $r^2=\rank(\shA)$, we have $\shA\simeq\O_X^{\oplus r^2}$ as 
$\O_X$-module.
The multiplication map $\shA\otimes\shA\ra \shA$ and the unit $\O_X\ra \shA$ 
induce   a $k$-algebra structure on $A=\Gamma(X,\shA)$. You easily check 
that $A\otimes \k(x)\simeq \shA(x)$ for each point $x\in X$.
So $A$ is a central simple $k$-algebra, which is trivial because $k$ is 
separably closed. Consequently $\shA=A\otimes\O_X$ is also trivial, 
contradiction.
\qed

\begin{remark}
It might easily happen that $X$ has trivial Picard group
\cite{Schroeer 1999a}. However, the preceding results ensures the existence 
of vector bundles of higher rank.
\end{remark}


\end{document}